\documentclass[a4paper]{article}

\usepackage[english]{babel}
\usepackage[utf8x]{inputenc}
\usepackage[T1]{fontenc}

\usepackage[a4paper,top=3cm,bottom=2cm,left=3cm,right=3cm,marginparwidth=1.75cm]{geometry}

\usepackage{amsmath}
\usepackage{graphicx}
\usepackage{hyperref}
\usepackage{xcolor}

\newtheorem{conj}{Conjecture}

\date{ }
\title{Notes on properties of binary strings which encode prime occurrences}
\author{Kajetan M\l ynarski\\\\
\small{Jagiellonian University, Krak\'{o}w, Poland}\\
\small{kajetan.mlynarski@uj.edu.pl}}

\begin{document}
\maketitle

\begin{abstract}

A binary string representation of prime occurrences is a sequence of bits, where $1$ entries encode positions of prime numbers. This is a convenient representation for analysis of prime distribution, since it allows for application of a broad range of existing string-analysis algorithms to problems in number theory. Binary strings of prime occurrences can be also generated with simple algorithms. Here we discuss three such algorithms and we demonstrate their applicability using the example of proving Goldbach's hypothesis for some limited sets of even numbers. This work formulates three open questions (conjectures) regarding the distribution of primes.

\end{abstract}


\section{Introduction}
\nocite{*}

Binary string of prime occurrences is a binary sequence where consecutive positions correspond to natural numbers, and non-zero entries occur in positions of prime numbers (e.g. $01101\ldots$ denotes first three primes: $2, 3, 5$). It is a very simple method of representing positions of primes in the set of natural numbers. This representation, while rarely applied (e.g. in visualizations of Ulam spirals \cite{Ulam}), is very clear and, importantly, can be easily computed with simple algorithms of low complexity. While we can not study the entire (infinite) string of prime occurrences, we can analyze regularities revealed by algorithms capable of generating this chain.
 
A particular type of algorithms capable of generating prime occurrence strings are incremental sieves. Incremental sieves allow for computation of $n+1$ prime, given the sequence of previous $n$ primes. An important property of incremental sieves is that they do not impose limits on the searched set of numbers, (as opposed to e.g. the Eratostenes sieve). General principles of incremental sieves are well known and used in construction of different algorithms \cite{Bengeloun, Pritchard, Sorenson}.

In this work we provide a synthesis of incremental sieves and binary strings of prime occurrences. We propose algorithms which generate binary strings of prime occurrences as  their output. This allows for simplification (through the reduction of computational complexity) of sieves. The goal of the proposed synthesis is however not fast computation of large primes, but development of a tool which could be useful in analysis of prime distribution. Developed algorithms could also find applications in other problems in number theory such as properties of sums of primes. Applications of sieve methods has a long tradition in number theory with applications to classical problems such as the Goldbach conjecture \cite{Chen}. We demonstrate that our method can also have applications to that problem, and use it to proving Goldbach conjecture for certain sets of even natural numbers.

\section{Algorithms}

\subsection{Basic algorithm - incremental sieve of Erathostenes}
\label{sec:basicAlgorithm}

This algorithm relies on two fundamental notions:

\begin{enumerate}
 \item Each number which is not a product of primes smaller than itself is a prime
 \item Each number which is a product of primes smaller than itself is composite
\end{enumerate}.

The algorithm proceeds as follows:

\begin{enumerate}
	\item Let $B$ be an infinite binary string initialized with of zeros only: $00000000 \ldots$
    $B$ is going to be transformed into a string of prime occurrences by replacing some of the $0$'s with $1$'s.
    
    \item In the line below $B$ - write $p_1$ - an infinite chain of occurrences of multiplicities of the smallest prime (i.e. $2$):\\
    $B: 00000000 \ldots$ \\
    $p_1: 01010101 \ldots$ \\
    
    $p_1$ is called the \emph{first sieve string}. 
    
    \item Flip $n-th$ entry of $B$ to $1$, where $n$ is the location of the first $1$ in the sieve string $p_1$: \\
    $B: 01000000 \ldots$ \\
    $p_1: 01010101 \ldots$ \\
    
    All nonzero $p_1$ entries following $n$ denote positions of composite numbers. First $0$ entry following $n$ indicates position of the consecutive prime (in our case $3$). The position of that entry is denoted by $m_1$.
    
    \item Flip $m_1-th$ entry of $B$ to $1$.
     
    \item Write consecutive $p_2$ sieve-string below. Nonzero entries in $p_2$ denote positions of multiplicities of the second smallest prime i.e. $3$, identified in the previous step: \\ 
     $B: 01000000 \ldots$ \\
     $p_1: 01010101 \ldots$ \\
     $p_2: 00100100 \ldots$ \\
     
     \item Position of the first column of zero entries ($m_2$) denotes position of the consecutive prime, analogous to $m_1$. 
   
\end{enumerate}

These steps can be further iterated to compute consecutive primes, until we obtain desired $k-th$ prime. In practical applications binary string $B$ does not have to be infinite. In such case, length of $B$ determines maximal range of searched primes.

\subsection{Algorithm 2 - quasi periodic}

The previous algorithm relied on sieve strings $p_1, p_2, \ldots$ to exclude multiplicities of consecutive primes. To compute $n-th$ prime, $n$ strings were necessary. Here, we note that these sieve strings can be superimposed, by applying a logical OR operation in each column. A superposition of sieve strings is called a \emph{cyclic sieve string}. An individual cyclic sieve $C_n$ string excludes multiplicities of all primes non larger than the $n-th$ one. For example, sieve strings $p_1, p_2$ can be superimposed in the following way: \\
$p_1: 01010101010\ldots$\\
$p_2: 00100100100\ldots$\\
$C_2: 011101011101...$\\
where the bottom string $C_2$ results from applying a logical OR operation to the strings above ($p_1, p_2$).

In such a way, we obtain single string $C_n$. $C_n$ consists of infinite cycles of a single sequence $S_n$ of length equal to primorial of $n$ ($p_n\#$). We denote $t-th$ occurrence of $S_n$ in $C_n$ as $S_{n,t}$. For consecutive $4$ primes $S_n$ has length equal to:
$2$\\
$2 \cdot 3 = 6$\\
$2 \cdot 3 \cdot 5 = 30$\\
$2 \cdot 3 \cdot 5 \cdot 7 = 210$\\

Each sequence $S_n$ selects all primes which are not multiplicities of $n$ first primes. Importantly, however $S_n$ does not select only primes, since even for $n=4$ it contains $1$ at a position corresponding to a composite number ($121$). Such composite numbers are included in the sequence, but are not eliminated. This is because $p_n\#$ grows faster than exponents of certain primes. We call these primes \emph{distracting numbers}.

For example: first occurrence of sequence $S_4$ in string $C_4$ indicates $121$ as a prime number. $121$ is composite (a square of $11$). $11$ in turn has been selected by $S_{4,1}$ ($11-th$ position of $S_4$ in its first occurrence in $C_4$ is equal to 1).

Existence of distracting numbers complicates the algorithm because their multiplicities have to be eliminated. This elimination can be implemented according to the following general rule:
all composite numbers belonging to interval $[1, p_n\#]$ are multiplicities
of $n$ first primes or of distracting numbers from that interval. Distracting numbers within the interval $[1, p_n\#]$ belong to the sub-interval $[z_n, \sqrt[]{p_n\#}]$, where $z_n$ is the $n-th$ prime.

We use general rule described above to modify the basic algorithm from section \ref{sec:basicAlgorithm} in two ways. The resulting variants of the algorithm highlight certain regularities in prime occurrences. The first variant relies 

\subsubsection{First variant}
\label{sec:firstMod}

\begin{enumerate}
	\item We initialize the algorithm with the string: $011010$, of length $p_2\#=6$. In this string positions of $1$ entries indicate primes (e.g. $2, 3$, which contribute to the primorial $p_2\#=6$) or other numbers which are not multiplicities of these primes. The smallest number which is not a multiplicity of preceding primes is itself a prime.
    
    \item The goal of the second step is the generation a consecutive string of length $p_3\#=30$. This is achieved by transforming the string from the first step according to the following rules:
    \begin{enumerate}
    	\item We replace $0$ in the first position with $1$s.
        \item We replace $1$ entries corresponding to factors of the primorial $p_2\#=6$  with $0$'s.
        \item Remaining positions remain unchanged. The result of the transformation is the string $100010$.
   		\item We append $z_{n}-1$ copies of the string generated in the step 2.3 to the end of the string generated in the 2.1 step, where $z_n$ is the $n-th$ prime. In this example, we append $4$ copies. The resulting string has the following form: $011010100010100010100010100010$. In this string, positions of the first two $1$ entries indicate primes $2,3$. The following $1$ entries indicate position of numbers, which are not multiplicities of $2, 3$. 
   		\item We replace $1$ entries corresponding to multiplicities of the 3rd prime (in a general case $n-th$ prime) with $0$. The $n-th$ prime is identified in step 1 above. We begin with the positions larger or equal to the square of the $n-th$ prime.
   \end{enumerate} 
   We obtain the following binary string of prime occurrences of length equal to $p_3\#$:\linebreak $011010100010100010100010100000$
   
   \item In the $n-th$ step we compute a $p_n\#$-long binary string of prime occurrences. This is done in a way analogous to the previous step:
   \begin{enumerate}
   		\item We transform the string obtained in the $n-1-th$ step by replacing $0$ in the first position with $1$.
        \item We replace $1$ entries corresponding to factors of the primorial $p_2\#=6$  with $0$'s.
        \item Remaining positions remain unchanged. 
        \item We append $m-1$ copies of the string generated in the step 3.3 to the end of the string generated in the step 3.1, where $m$ is the $n-th$ prime.  
        \item We replace $1$ entries corresponding to multiplicities of the $n-th$ prime (in a general case $n-th$ prime) with $0$. The $n-th$ prime is identified in step 1 above. We begin with the positions larger or equal to the square of the $n-th$ prime. We repeat this procedure for all numbers indicated by $1$ entries in the interval $[z_n, \sqrt[]{p_n\#}]$ i.e. all distracting numbers. It is important to note, that for each distracting number entries corresponding only to either its exponents or products of a that distracting number and any prime from the $[z_n, \sqrt[]{p_n\#}]$ interval have to be modified.
   \end{enumerate}  
\end{enumerate}

The above steps can be iterated until the desired prime has been generated.

\subsubsection{Second variant}

We begin by initializing the algorithm with the $p_n\#$ long binary string of prime occurrences. This "seed" string can be obtained with the algorithms described above. $n$ is primarily bounded by accessible computational resources. We then transform the seed string according to the rule described in the third step of the first modification (\ref{sec:firstMod}). The string resulting from that transformation is appended to the seed string $q$ times. The exact value of $q$ is specified by how large prime we want to obtain. We then replace all $1$ entries on positions corresponding to multiplicities of all primes from the interval $[z_n, \sqrt[]{q p_n\#}]$ with $0$'s.  


\section{Representations}

A binary representation of a natural number $n$ is a binary string of length n, where positions of $1$’s and $0$’s denote natural numbers smaller or equal to $n$ which satisfy some pre-defined properties. In particular, $1$'s and $0$'s can correspond to primes and composite numbers smaller than n respectively. For example, $011010$ represents $6$ in such representation. An example of a different representation of $6$ is a binary string $010101$ where 1's represent positions of even numbers smaller or equal to $6$. These representations provide a useful tool for analysis of number properties by displaying certain aspects of their structure.

A particularly interesting application of binary representations of natural numbers is analysis of properties of natural numbers  (e.g. their sums and products) without the use of arithmetic.

Let us write down a binary representation of 24 and mark all multiplicities of 2: \newline
$010101010101010101010101$ \linebreak
Above, let us write down the same string in the reverse order, from 24 to 1: \newline
 $10101\textcolor{red}{0}101010101010101010$\linebreak
 $0101\textcolor{red}{0}1010101010101010101$\linebreak
The two to odd numbers 5 and 19 (numbers in the top string) are marked in red. The sum of 5 and 19 is obviously 24. Each 0 entry in the bottom string has a corresponding 0 in the top string. The position of this corresponding 0 in the top row (counted from the right side) when added to the position of the 0 entry in the bottom row yields 24. This relationship holds for all binary string representations of even numbers.

This implies that for every odd number $P$ smaller than an even number $M$ (such that $|P-M| > 2$) there exists an odd number $P'$ such that: $P + P' = M$ (except for $N=4$) is a sum of all numbers denoted by positions of zeros in the bottom string and zeros closest to them. Therefore every even number larger than $4$ is a sum of two odd numbers and is not a sum of an odd number and $2$. 

This extremely simple example demonstrates an application of binary string representations to analyze properties of natural numbers. 

Let us consider a string consisting of a single cycle of length equal to $p_3\#=2 \cdot 3 \cdot 5 = 30$ with $1$ entries at positions corresponding to prime numbers. This is a binary string representation of the number $30$ with respect to positions of primes: $011010100010100010100010000010$.

All these prime numbers are determined by elimination of multiplicities of $2, 3, 5$ (the issue of distracting numbers does not occur here).

As in the example above, let us write the inverted version of the same string in the row above:\linebreak
$010000010001010001010001010110$\linebreak
$011010100010100010100010000010$\linebreak
Analogously to the previous example we will exclude composite numbers. Because excluded composite numbers are equal to multiplicities of $2, 3, 5$, therefore $30$ can not be a sum of $2, 3, 5$ and another prime number. This is because $2, 3, 5$ are always added to its own multiplicities, which are by definition composite. Each remaining prime, is therefore paired with another remaining prime, such that their sum is equal to $30$. This is analogous to the previous example, in which we demonstrated that no even number with the exception of $4$ can not be a sum of $2$ and an odd number.

It is now easy to observe that no even number is a sum of its own prime factor and a prime number. The only exception are numbers of the form $2 \times M$, where $M$ is a prime number.

We can now formulate more general criteria under which prime numbers smaller than $p_n\#$ do not sum to $p_n\#$ when added to another prime. This true in the following cases:
\begin{enumerate}
    
    \item Prime factors of $p_n\#$. They sum to $p_n\#$ when added to their own multiplicities, which are by definition composite.
    
    \item Prime numbers smaller than $p_n\#$, which sum up to $p_n\#$ when added to multiplicities of distracting numbers.
    
\end{enumerate}
Every prime number which does not meet the above criteria sums up to $p_n\#$ when added to another prime number. 

In case of numbers greater than $p_3\#$ there exist distracting numbers, which belong to the interval $[z_n, \sqrt{p_n\#}]$, where $z_n$ is the larges prime contributing to $p_n\#$. 

Let us observe the following:
\begin{itemize}
    \item $\displaystyle{\lim_{n \to \infty}} \frac{\sqrt{p_n\#}}{p_n\#} = 0$
    
    \item The proportion of the distracting numbers in the set of primes smaller than $p_n\#$ is decreasing with increasing $n$. Moreover, it has been verified numerically \cite{Silva} that for $n \leq 15$, $p_n\#$ is a sum of two prime numbers.
    
\end{itemize}
Therefore for each number of the form $p_n\#$ exists at least one pair of prime numbers $p, p'$ such that: $p + p' = p_n\#$, i.e. every number of the form $p_n\#$ is a sum of two prime numbers.

Let us consider numbers of the following form: $m p_n\#$. These numbers are even and each of them is a sum of prime numbers $p, p'$. Individually, numbers $p, p'$ satisfy the following properties:

\begin{enumerate}
    \item They are not prime factors of $m p_m\#$.
    \item The do not sum up to $m p_m\#$ when added to multiplicities of distracting numbers.
\end{enumerate}

In the case of numbers of the form $m p_n\#$ distracting numbers belong to the interval $[z_n, \sqrt{m p_n\#}]$. Similarly as in the case of the numbers of the form $p_n\#$, $\displaystyle{\lim_{n \to \infty}} \frac{\sqrt{m p_n\#}}{mp_n\#} = 0$. Because of that fact, for each number of the form $m p_n\#$ there exists at least one pair of prime numbers $p, p'$ such that: $p + p' = m p_n\#$. Therefore every number of the form $m p_n\#$ is a sum of two primes.

The smallest number of the form $p_n\#$ is $p_1\# = 2$. All even numbers are of the form $p_1\#$. Therefore for each even number there exits at least one pair of prime numbers $p, p'$ such that: $p + p' = m p_n\#$. Each even number is a sum of two primes.

\section{Remarks and conjectures about distribution of primes}

Let us consider the following two binary strings:\\
$001101010001010001010001$ \linebreak
$011010100010100010100010$ \linebreak 
In the bottom string of length equal to $24$ positions of all primes smaller than $24$ were denoted with $1$. In the top row positions of all numbers 

If any of the second summands is a prime number, then the given number (in that case $24$) is a sum of two primes. An even number therefore is not a sum of two primes only if all second summands are not prime.

The string of second summands is generated by shifting the string of primes by one position to the right (adding a $0$ at the beginning and deleting the last position). In contrast to the bottom string (a binary string of prime occurrences), positions in the string of second summands are indexed from the right. Indices of second summands are therefore odd and are distributed in the same distance as prime numbers. If all second summands were composite, all distances between their positions in the string would have to be equal to multiplicities of prime numbers smaller than $\sqrt{n}$, where $n$ is a given even number (in our case $n=24$). This, however is impossible, since distances between second summands are equal to distances between prime numbers. This implies that at least one second summand must be a prime. If all second summands were composite, then one of the primes smaller than $n$ could not be a prime. 

Let $d_{\max}$ be the largest distracting number. 

\begin{conj}
For each $n$ in the interval $[p_n\# - d_{max}, p_n\# + d_{max}]$ there exists at least one prime.
\end{conj}

\begin{conj}
At least one prime from the interval $[p_n\# - d_{max}, p_n\# + d_{max}]$ is a sum or a difference $p_n\#$ and a prime from the interval $[1, d_{\max}]$.
\end{conj}

\begin{conj}
The prime closest to $p_n\#$ is equal to a sum of $p_n\#$ and a number from the interval $[1, d_{\max}]$
\end{conj}

\pagebreak

\bibliographystyle{alpha}
\bibliography{main}

\end{document}